\def\draft{n}
\documentclass{amsart}
\usepackage[headings]{fullpage}
\usepackage{amssymb} 

\theoremstyle{plain}

\newtheorem{theorem}{Theorem}
\newtheorem{proposition}{Proposition}[section]
\newtheorem{lemma}[proposition]{Lemma}

\theoremstyle{definition}

\theoremstyle{remark}

\newtheorem{remark}[proposition]{Remark}

\def\printname#1{
        \if\draft y
                \smash{\makebox[0pt]{\hspace{-0.5in}
                        \raisebox{8pt}{\tt\tiny #1}}}
        \fi
}

\newlength{\standardunitlength}
\catcode`\@=11
\long\def\@makecaption#1#2{%
     \vskip 10pt

\setbox\@tempboxa\hbox{
       \small\sf{\bfcaptionfont #1. }\ignorespaces #2}%
     \ifdim \wd\@tempboxa >\captionwidth {%
         \rightskip=\@captionmargin\leftskip=\@captionmargin
         \unhbox\@tempboxa\par}%
       \else
         \hbox to\hsize{\hfil\box\@tempboxa\hfil}%
     \fi}
\font\bfcaptionfont=cmssbx10 scaled \magstephalf
\newdimen\@captionmargin\@captionmargin=2\parindent
\newdimen\captionwidth\captionwidth=\hsize
\catcode`\@=12

\newcommand{\tr}{\operatorname{tr}}

\def\lbl#1{\label{#1}\printname{#1}}

\def\BN{\mathbb N}

\def\P{\mathcal P}
\def\R{\mathcal R}

\def\La{\Lambda}

\def\d{\delta}

\renewcommand{\det}{\mathrm{det}}
\def\inv{\mathrm{inv}}
\def\Ferm{\mathrm{Ferm}}
\def\Bos{\mathrm{Bos}}

\begin{document}


\title[The quantum MacMahon Master Theorem]{
The quantum MacMahon Master Theorem}

\author{Stavros Garoufalidis}
\address{School of Mathematics \\
         Georgia Institute of Technology \\
         Atlanta, GA 30332-0160, USA \\ 
         {\tt http://www.math.gatech} \newline {\tt .edu/$\sim$stavros } }
\email{stavros@math.gatech.edu}
\author{Thang TQ L\^e}
\address{Department of Mathematics \\
         Georgia Institute of Technology \\
         Atlanta, GA 30332-0160, USA}
\email{letu@math.gatech.edu}
\author{Doron Zeilberger}
\address{Department of Mathematics \\
         Rutgers University \\ 
         110 Frelinghuysen Rd \\
         Piscataway, NJ 08854-8019, USA}
\email{zeilberg@math.rutgers.edu}

\thanks{The authors were supported in part by NSF.
\newline
1991 {\em Mathematics Classification.} Primary 57N10, 05A30. Secondary 57M25.
\newline
{\em Key words and phrases: Boson-fermion correspondence, MacMahon's
Master Theorem, $q$-difference operators} 
}

\date{
This edition: May 26, 2005. 
\hspace{0.5cm} First edition: March 24, 2003.}

\begin{abstract}
We state and prove a quantum-generalization of MacMahon's celebrated Master 
Theorem, and relate it to a quantum-generalization of the boson-fermion
correspondence of Physics.
\end{abstract}

\maketitle



\section{Introduction}
\lbl{sec.intro}
\subsection{MacMahon's Master Theorem}
\lbl{sub.mm}

In this paper we state and prove a quantum-generalization 
of MacMahon's celebrated Master Theorem, conjectured by the first two 
authors. Our result was motivated by quantum topology. 
In addition to its potential importance in
knot theory and quantum topology (explained in brief in the last section), 
this paper answers 
George Andrews's long-standing open problem \cite{A} of finding a
{\it natural} q-analog of MacMahon's Master Theorem.

Let us recall the original form of MacMahon's Master Theorem and some of
its modern interpretations.

Consider a square matrix $A=(a_{ij})$ of size $r$ with entries in some 
commutative ring. For $1 \leq i \leq r$, let
$X_i:=\sum_{j=1}^{r} a_{ij} x_j  $, (where $x_i$'s are commuting variables)
and for any vector $(m_1, \dots, m_r)$ of non-negative integers
let $G(m_1, \dots, m_r)$ be 
the coefficient of $x_1^{m_1} x_2^{m_r} \dots x_r^{m_r}$ in
$\prod_{i=1}^r X_{i}^{m_i} $. 
{\em MacMahon's Master Theorem} is the following identity (see \cite{MM}):
\begin{equation}
\lbl{eq.mm}
\sum_{m_1, m_2, \dots, m_r =0}^\infty
G(m_1, \dots, m_r) 
= 1/\det(I-A) .
\end{equation}

There are several equivalent reformulations of MacMahon's Master Theorem;
see for example \cite{FZ} and references therein. Let us mention one, of
importance to physics.

Given a matrix $A=(a_{ij})$ of size $r$ with commuting entries
which lie in a ring $\R$, and a nonnegative integer $n$,
we can consider its {\em symmetric} and  {\em exterior}
powers $S^n(A)$ and $\La^n(A)$, and their traces $\tr S^n(A)$
and $\tr \La^n(A)$ respectively.
Since
\begin{eqnarray*}
\tr S^n(A) &=& \sum_{m_1+ \dots m_r=n} G(m_1,\dots,m_r) \\
\det(I-t A) &=& \sum_{n=0}^\infty (-1)^n \tr \La^n(A) t^n,
\end{eqnarray*}
the following identity
\begin{equation}
\lbl{eq.BF}
\frac{1}{\sum_{n=0}^\infty (-1)^n \tr \La^n(A) t^n}
=\sum_{n=0}^\infty \tr S^n(A) t^n
\end{equation}
in $\R[[t]]$ is equivalent to \eqref{eq.mm}. In
Physics \eqref{eq.BF} is called the {\em boson-fermion} correspondence,
where bosons (resp. fermions) are commuting (resp. skew-commuting) 
particles corresponding to symmetric (resp. exterior) 
powers.

\subsection{Quantum algebra, right-quantum matrices and quantum determinants}
\lbl{sub.qqq}

In $r$-dimensional 
{\it quantum algebra} we have $r$ indeterminate variables
$x_i$ ( $ 1 \leq i \leq r$), satisfying the
commutation relations $x_j x_i=qx_i x_j$ for all $1 \leq i<j \leq r$.
We also consider matrices $A=(a_{ij})$ of $r^2$ indeterminates $a_{ij}$, 
$1 \leq i, j \leq r$, that commute with the $x_i$'s
and such that for any $2$ by $2$  minor of
$(a_{ij})$, consisting of rows $i$ and $i'$,
and columns $j$ and $j'$ (where $1 \leq i < i' \leq r$,
and $1 \leq j < j' \leq r$), writing
$a:=a_{ij}, b:=a_{ij'}, c:=a_{i'j}, d:=a_{i'j'}$,
we have the {\em commutation relations}:
\begin{eqnarray}
\lbl{eq.commutation1}
ca &= & qac, \quad \text{($q$-commutation of the entries in a column)} \\
\lbl{eq.commutation2}
 db &=& qbd, \quad \text{($q$-commutation of the entries in a column)} \\
\lbl{eq.commutation3}
ad &=& da+q^{-1}cb-q bc \qquad \text{(cross commutation relation)}.
\end{eqnarray}
We will call such matrices $A$ {\em right-quantum matrices}.

The  {\em quantum determinant}, 
(first introduced in \cite{FRT})  of any (not-necessarily right-quantum)
$r$ by $r$ matrix $B=(b_{ij})$ may be defined by

$$
\det_q(B):=
\sum_{\pi \in S_r} (-q)^{-\inv(\pi)} b_{\pi_1 1} b_{\pi_2 2}
\cdots b_{\pi_r r} ,
$$
where the sum ranges over the set of permutations, $S_r$, of 
$\{1, \dots , r\}$,
and for any of its members, $\pi$, 
$\inv(\pi)$ denotes  the number of pairs $1 \leq i<j \leq r$
for which $\pi_i>\pi_j$.

\subsection{A $q$-version of MacMahon's Master Theorem}
\lbl{sub.statement}

We are now ready to state our quantum version of Macmahon's Master
Theorem.

\begin{theorem}
\lbl{thm.1}
{\bf (Quantum MacMahon Master Theorem)}
Fix a right-quantum matrix $A$ of size $r$. For $1 \leq i \leq r$, let
$X_i:=\sum_{j=1}^{r} a_{ij} x_j  $,
and for any vector $(m_1, \dots, m_r)$ of non-negative integers
let $G(m_1, \dots, m_r)$ be 
the coefficient of $x_1^{m_1} x_2^{m_r} \dots x_r^{m_r}$ in
$\prod_{i=1}^r X_{i}^{m_i} $.
Let
$$
\Ferm(A)=\sum_{J \subset \{1, \dots, r\}} (-1)^{|J|} \det_q(A_J)
$$
where the summation is over the set of {\em all} subsets $J$
of $\{1,\dots,
r\}$, and $A_J$ is the $J$ by $J$ submatrix of $A$, and
$$
\Bos(A)=\sum_{m_1,\dots, m_r=0}^\infty G(m_1,\dots,m_r).
$$
Then 
$$
\Bos(A)=1/\Ferm(A) .
$$
\end{theorem}

When we specialize to $q=1$, Theorem \ref{thm.1} recovers Equation
\eqref{eq.BF}, which explains why our result is a $q$-version of the
MacMahon Master Theorem. 
For a motivation of Theorem \ref{thm.1}, see Section \ref{sec.bf}.

The above result is not only interesting from the combinatorial point of
view, but it is also a key ingredient
in a {\em finite noncommutative formula} for the colored Jones function of 
a knot. This will be explained in a subsequent publication, \cite{HL}.

\subsection{Computer code}

The results of the paper have been verified by computer code, written 
by the third author. {\tt Maple} programs {\tt QuantumMACMAHON} and {\tt
qMM} are available at: {\tt http://www.math.rutgers.edu/}
\newline {\tt $\sim$zeilberg/}.
The former proves rigorously Theorem \ref{thm.1} for any fixed $r$.

\subsection{Acknowledgement}
The authors wish to thank the anonymous referee who pointed out an
error in an earlier version of the paper, and Martin Loebl for enlightening
conversations.

\section{Proof}

\subsection{Some lemmas on operators}
\lbl{sub.some}

The proof will make crucial use of a {\em calculus of difference operators},
developed by the third author in \cite{Z1}. This calculus of difference
operators predates the more advanced calculus of holonomic functions,
developed by the third author in \cite{Z2}.

Difference operators act on
{\em discrete functions} $F$, that is functions whose domain is $\BN^r$.
For example, consider the {\it shift-operators} $M_i$ and the {\em
multiplication} operator $Q_i$
which act on a discrete function $F(m_1, \dots, m_r)$ by
\begin{eqnarray*}
(M_iF)(m_1, \dots, m_r) &:=&
F(m_1, \dots, m_{i-1}, m_i+1, m_{i+1}, \dots, m_r) \\
(Q_iF)(m_1, \dots, m_r) &:=& q^{m_i} F(m_1, \dots, m_r).
\end{eqnarray*}
It is easily seen that
$$
M_i Q_i = q Q_i M_i.
$$
Abbreviating $Q_i$ by $q^{m_i}$, we obtain that:
\begin{equation}
\lbl{eq.EQ}
M_i q^{m_i}=
q^{m_i+1} M_i  \qquad M_i q^{m_j}=q^{m_j} M_i \quad \text{for} \quad
i \neq j.
\end{equation}

Another example is the operator $\hat{x}_i$ which left multiplies $F$ by $x_i$.
Notice that $\hat{x_j} \hat{x}_i = q \hat{x}_i \hat{x}_j$ for $j>i$.
In the proof below, we will denote $\hat{x}_i$ by $x_i$. In that case,
the identity $x_j x_i=q x_i x_j$ for $j>i$ holds in the
quantum algebra, as well as in the algebra of operators.

Before embarking on the proof,
we need the following readily-verified lemmas.

\begin{lemma}
\lbl{lem.1} (commuting $X_i$ with $X_j$)
For $1 \leq i < j \leq r$, $X_jX_i=qX_iX_j$.
\end{lemma}

\begin{lemma}
\lbl{lem.2} (commuting $x_i$ with $X_j$)
For each of the $a_{ij}$, define the
operator $Q_{ij}$ acting on expressions $P$ involving $a_{ij}$
by $Q_{ij} P(a_{ij}):=P(qa_{ij})$.
Then, for any $1\leq i,j \leq r$, and integer $m_i$ and any expression $F$
$$
x_i^{-m_i}X_j F=
[(Q_{j1}^{-1}Q_{j2}^{-1} \cdots Q_{j,i-1}^{-1}  Q_{j,i+1} \cdots
Q_{jr})^{m_i}X_j] x_i^{-m_i} F .
$$
\end{lemma}

\begin{lemma}
\lbl{lem.3}
(Column expansion with respect to the last column):
Given an $r$ by $r$ matrix $(a_{ij})$ (not necessarily quantum) let
$A_i$ be the minor of the entry $a_{ir}$, i.e. the
$r-1$ by $r-1$ matrix obtained by deleting the $i^{th}$ row
and  $r^{th}$ column. Then
$$
\det_q(A)=\sum_{i=1}^{r} (-q)^{i-r} \left ( \det_q{A_i} \right ) a_{ir}
 .
$$
\end{lemma}


\begin{lemma}
\lbl{lem.3a}
If $A$ is a matrix that satisfies Equation
\eqref{eq.commutation3} and $A'$ denotes a matrix obtained
by interchanging the $i$ and $j$ columns columns of $A$, 
then $\det_q(A')=(-q)^{-\inv(ij)} \det_q(A)$.
\end{lemma}

\begin{proof}
Suppose first that we interchange two adjacent colums $i$ and $j:=i+1$.
Consider the involution 
of $S_r$ that sends a permutation $\pi$ to $\pi'=\pi(ij)$. Given
$\pi \in S_r$, let $(A;\pi)=(-1)^{-\inv(\pi)} a_{\pi_1 1} \dots a_{\pi_r r}$
denote the contribution of $\pi$ in $\det_q(A)$. Then, $\det_q(A)=\sum_{\pi}
(A;\pi)$. Equation \eqref{eq.commutation3} implies that
$$
(A;\pi)+(A;\pi')=(-q)((A';\pi) + (A';\pi')).
$$
Summing over all permutations proves the result when $j=i+1$.

Observe that when $j=i+1$, the matrix $A'$ is no longer right-quantum
since it does not satisfy \eqref{eq.commutation3}.
However, the proof used only the fact that \eqref{eq.commutation3}
holds for the $i$ and $i+1$ columns of $A$.

Thus, the proof can be iterated $\inv(ij)$ times 
to commute the $i$ and $j>i$ columns of $A$. The result follows.
\end{proof}

\begin{lemma}
\lbl{lem.4} (Equal columns imply that $\det_q$ vanishes):
Let $A$ be a right-quantum matrix. 
In the notation of Lemma \ref{lem.3}, for all $j \neq r$,
$$
\sum_{i=1}^{r} (-q)^{i-r} \left ( \det_q{A_i} \right ) a_{ij}=0 .
$$
\end{lemma}

\begin{proof}
If $j=r-1$, it is easy to see that $q$-commutation along
the entries in every column of $A$ imply that the sum vanishes.

If $j < r-1$, use Lemma \ref{lem.3a} to reduce it to the 
case of $j=r-1$. 
\end{proof}

\begin{remark}
\lbl{rem.laplace}
One can give an alternative proof of Lemmas \ref{lem.3a} and \ref{lem.4}
from the trivial $2$ by $2$ case
and, by induction
using the $q$-Laplace expansion of a $q$-determinant 
that is completely analogous
to the classical case.
\end{remark}

\subsection{Proof of Theorem \ref{thm.1}}
\lbl{sub.proof}

The proof is a quantum-adaptation
of the ``operator-elimination" proof of MacMahon's Master Theorem
given in \cite{Z1}. Fix a right-quantum matrix $A$.

Observe that
$G(m_1, \dots, m_r)$ is the coefficient of $x_1^{0} \dots x_r^{0}$ in
$$
H(m_1, \dots, m_r; x_1, \dots, x_r):=
x_r^{-m_r} \cdots x_2^{-m_2} x_1^{-m_1}\prod_{i=1}^r X_{i}^{m_i} .
$$
We will think of $H$ as a {\em discrete function}, that is as a function
of $(m_1,\dots,m_r) \in \BN^r$. $H$ takes values in the ring of
noncommutative Laurrent polynomials in the $x_i$s, with coefficients in
the ring generated by the entries of $A$,
modulo the ideal given by \eqref{eq.commutation1}-\eqref{eq.commutation3}.

Let's see how the shift operators $M_i$ acts on $H$. By definition,
$$
M_i H(m_1, \dots, m_r; x_1, \dots, x_r)=
$$
$$
x_r^{-m_r} \cdots x_{i+1}^{-m_{i+1}} x_{i}^{-m_{i}-1}
x_{i-1}^{-m_{i-1}} \cdots x_1^{-m_1}
X_1^{m_1} \cdots X_{i-1}^{m_{i-1}} X_{i}^{m_{i}+1}
X_{i+1}^{m_{i+1}} \cdots X_{r}^{m_{r}}  .
$$
By moving ${x_i}^{-1}$ to the front and 
$X_i$ in front of ${X_1}^{m_1}$, and using Lemma \ref{lem.1}
and $x_jx_i=qx_ix_j$, we have
$$
M_i H(m_1, \dots, m_r;x_1, \dots, x_r)=
q^{m_r+m_{r-1}+ \dots + m_{i+1}- m_1 - m_2 - \dots -m_{i-1}}
x_i^{-1} [ x_r^{-m_r} \cdots x_1^{-m_1}  X_i]
X_1^{m_1} \cdots X_{r}^{m_{r}}  .
$$
By moving $X_i$ next to $x_i^{-1}$ and using Lemma \ref{lem.2}
this equals to:
$$
q^{m_r+m_{r-1}+ \dots + m_{i+1}- m_1 - m_2 - \dots -m_{i-1}}x_i^{-1} \cdot
$$
$$
[(Q_{i2} \cdots Q_{ir})^{m_1}
(Q_{i1}^{-1} Q_{i3} \cdots Q_{ir})^{m_2}
(Q_{i1}^{-1} Q_{i2}^{-1} Q_{i4} \cdots Q_{ir})^{m_3} \cdots
(Q_{i1}^{-1} Q_{i2}^{-1}  \cdots Q_{i,r-1}^{-1})^{m_r} X_i ] \cdot
$$
$$
x_r^{-m_r} \cdots x_1^{-m_1}{X_1}^{m_1} \dots {X_r}^{m_r}  ,
$$
which is equal to
$$
q^{m_r+m_{r-1}+ \dots + m_{i+1}- m_1 - m_2 - \dots -m_{i-1}}{x_i}^{-1}
\cdot
$$
$$
\left ( q^{-m_2-m_3- \dots -m_r}a_{i1}x_1+
q^{m_1-m_3- \dots -m_r}a_{i2}x_2+ \dots+
q^{m_1+m_2+\dots +m_{r-1} }a_{ir}x_r  \right )
H(m_1, \dots, m_r; x_1, \dots, x_r).
$$

Multiplying out and rearranging, we get that the discrete function
$H(m_1, \dots, m_r; x_1, \dots, x_r)$  is annihilated by the
$r$ operators ($i=1,2, \dots, r$)
$$
\P_i:=
\sum_{j=1}^{i-1} -q^{-m_j-2m_{j+1}-\dots -2m_{i-1}-m_i}a_{ij}x_j+
(M_i-a_{ii})x_i+
\sum_{j=i+1}^{r} -q^{m_i+2m_{i+1}+\dots +2m_{j-1}+m_j}a_{ij}x_j .
$$

Now comes a nice surprise. Let us define $b_{ij}$ to be
the coefficient of $x_j$ in $\P_i$. For example, for $r=3$ we have:
$$
B= \left(
\begin{matrix}
M_1-a_{11} & -q^{m_1+m_2} a_{12} & -q^{m_1+2m_2+m_3} a_{13} \\
 -q^{-m_1-m_2} a_{21} & M_2 -a_{22} & -q^{m_2+m_3} a_{23} \\
-q^{-m_1-2m_2-m_3} a_{31} & -q^{-m_2-m_3} a_{32} & M_3 -a_{33} 
\end{matrix} \right).
$$

\begin{lemma}
\lbl{lem.5}
$B$ is a right-quantum matrix.
\end{lemma}

\begin{proof}
It is easy to see that the entries in each column of $B$ $q$-commute.
To prove Equation \eqref{eq.commutation3}, consider the 
following cases for a $2$ by $2$ submatrix $C$ of $B$: $C$ contains
two, (resp. one, resp.  no) diagonal entries of $B$, and prove it case
by case, using the fact that the operators $M_i$ and $q^{m_j}$ commute with the
$a_{ij}$, and satisfy the commutation relations \eqref{eq.EQ}.
\end{proof}

Now we eliminate $x_1, x_2, \dots , x_{r-1}$ by left-multiplying
$\P_i$ by the minor of $b_{ir}$ in 
$B=(b_{ij})$
times $(-q)^{i-r}$, for each $i=1,2, \dots, r$, and adding them all up.
Since $B$ is right-quantum (by Lemma \ref{lem.5}), 
Lemma \ref{lem.4} implies that 
the coefficients of $x_1, \dots, x_{r-1}$ all vanish,
and $\det_q(B) x_r H=0$. After left multiplying by $x_r^{-1}$ which
commutes with the entries in $B$, we obtain that
$$
\det_q(B) H(m_1, \dots, m_r;x_1,\dots,x_r)=0.
$$
Since the entries of $B$ do not contain $x_i$'s, it follows that
$\det_q(B)$ annihilates
every coefficient of $H$, in particular its constant
term. Taking the constant term yields
$$
\det_q(B) G(m_1, \dots, m_r)=0 .
$$
Here comes the next surprise.
\begin{lemma}
\lbl{lem.6}
\rm{(a)}
We have:
$$
\det_q(B)=\sum_{J \subset \{1,\dots,r\}} (-1)^{|J|} \det_q(A_J) M_{\bar J}
$$
where $\bar J=\{1,\dots,r\}-J$ and $M_J=\prod_{j \in J} M_j$. \newline
\rm{(b)} In particular,
$$
\det_q(B)|_{M_1=\dots=M_r=1}
=\Ferm(A).
$$
\end{lemma}

\begin{proof}
Let us expand $\det_q(B)$ as a sum over permutations $\pi \in S_r$.
We have:
\begin{eqnarray*}
\det_q(B) &=& \sum_{\pi \in S_r} (-q)^{-\inv(\pi)} b_{\pi_1 1} b_{\pi_2 2}
\cdots b_{\pi_r r} \\
&=& \sum_{\pi \in S_r} \prod_{i=1}^r (-q)^{-\inv(\pi,i)} b_{\pi_i i} 
\end{eqnarray*}
where
$\inv(\pi,i)$ is the number of $j >i$ such that $\pi_i > \pi_j$. 
Now, $b_{ij}=\d_{ij}M_i-q_{ij} a_{ij}$,
where $q_{ij}$ is a monomial in the variables $q^{m_k}$, 
and $\prod_i q_{\pi_i i}=1$. Moreover, if $\pi_i=i$, then for each $j$
with $i < j \neq \pi_j$, the exponent of $q^{m_i}$ in
$q_{ij}$ is $2$ if $\pi_j < i$ and $0$ if $\pi_j > i$.

Since $\prod_i q_{\pi_i i}=1$, we can move the monomials $q_{ij}$ 
in the left of $\prod_i (-q)^{-\inv(\pi,i)}
b_{\pi_i i}$, and then cancel them. The monomials commute
with all entries of the matrix $b_{ij}$, {\em except} with the diagonal ones.
Commuting $q^{2 m_i}$ with $b_{ii}=\d_{\pi_i i} M_{ii} - 
q_{\pi_i i} a_{\pi_i i}$ gives:
$b_{ii} q^{2 m_i}=q^{2 m_i} (\d_{\pi_i i} q^2 M_{ii} - 
q_{\pi_i i} a_{\pi_i i})$. In other words, it replaces $M_i$ by $q^2 M_i$.
Thus, we have:
\begin{eqnarray*}
\det_q(B) &=&
\sum_{\pi \in S_r} \prod_{i=1}^r (-q)^{-\inv(\pi,i)} (\d_{\pi_i i}
q^{2 \inv(\pi,i)} M_i-
a_{\pi_i i}) \\
&=& 
\sum_{\pi \in S_r} \sum_{J \subset \{1, \dots, r\}}
\prod_{i \in J} (-q)^{-\inv(\pi,i)} \d_{\pi_i i}
q^{2 \inv(\pi,i)} M_i \prod_{i \not\in J} (-q)^{-\inv(\pi,i)} (-a_{\pi_i i})
\end{eqnarray*}
Now, rearrange the summation. Observe that
every permutation $\pi$ of $\{1,\dots,r\}$ gives rise to a permutation
$\pi'$ on the set $\{1,\dots,r\}-\text{Fix}(\pi)$, where $\text{Fix}(\pi)$
is the fixed point set of $\pi$. Moreover, $\inv(\pi',i)=\inv(\pi,i)-
|\{j \in J: j>i\}|$. Using this, part (a) follows. Part (b) follows from
part (a) and the definition of $\Ferm(A)$.
\end{proof}

Hence
$$
\sum_{J \subset \{1,\dots,r\}} (-1)^{|J|} \det_q(A_J) M_{\bar J}
 G(m_1, \dots, m_r)=0 .
$$
Summing over $\BN^r$, we get:
$$
\sum_{m_1,\dots,m_r=0}^\infty
\sum_{J \subset \{1,\dots,r\}} (-1)^{|J|} \det_q(A_J) M_{\bar J}
G(m_1, \dots, m_r)=0 .
$$
For a subset $J=\{k_1,\dots,k_j\}$ of $\{1,\dots, r\}$, we denote by
$G_J(m_{k_1},\dots,m_{k_j})$ the evaluation $G(m_1,\dots,m_r)$
at $m_i=0$ for all $i \not\in J$, and we define
$$
S_J=\sum_{m_{k_1}, \dots, m_{k_j}=0}^\infty G(m_1,\dots,m_r).
$$
Using telescoping cancellation, the inclusion-exclusion principle,
and Lemma \ref{lem.6}(b), the above equation becomes
$$
\sum_{J \subset \{1,\dots,r\}}
(-1)^{|J|} \text{Ferm}(A_J) S_J=0.
$$
Using induction (with respect to $r$), together with $S_{\emptyset}=1$,
we obtain that
$\Ferm(A) S_{\{1,\dots,r\}}=1$. This concludes the proof of the theorem.
\qed

\section{Some remarks on the boson-fermion correspondence}
\lbl{sec.bf}

Let us give some motivation for Theorem \ref{thm.1} from the
point of view of quantum topology.

For a reference on quantum space
and quantum algebra, see \cite[Chapter IV]{Ka} and \cite{M}.

Recall that a vector (column or row) of $r$ indeterminate entries 
$x_1,\dots,x_r$ lies in  $r$-dimensional {\em quantum space} $A^{r|0}$ iff 
its entries satisfy
$$ 
x_j x_i=qx_i x_j
$$ 
for all $1 \leq i<j \leq r$. 

Recall that a {\em right} (resp. {\em left}) {\em endomorphism} 
of $A^{r|0}$ is a matrix
$A=(a_{ij})$ of size $r$ whose entries commute with the coordinates
$x_i$ of a vector $x=(x_1,\dots,x_r)^T \in A^{r|0}$ and in addition,
$Ax$ (resp. $x^T A$) lie in $A^{r|0}$. Recall also that an endomorphism
of $A^{r|0}$ is one that is right and left endomorphism.

It is easy to see (eg. in \cite[Thm. IV.3.1]{Ka})
that $A$ is a right-quantum (i.e., a right-endomorphism) 
iff  for every $2$ by $2$ submatrix $\begin{pmatrix}
a & b \\ c & d \end{pmatrix}$ of $A$ we have:
$$
ca=qac, \qquad db=qbd,  \qquad
\quad ad=da+q^{-1}cb-q bc.
$$
Moreover, $A$ is left-quantum iff for every $2$ by $2$ submatrix of $A$
(as above) we have:
$$
ba=qab, \qquad dc=qcd, \qquad  ad=da+q^{-1}bc-q cb.
$$
Finally, $A$ is quantum iff for every $2$ by $2$ submatrix of $A$ (as above)
we have:
\begin{equation}
\lbl{eq.commutation}
ba=qab, \quad ca=qac, \quad db=qbd, \quad dc=qcd, \quad
cb=bc, \quad ad=da+q^{-1}cb-q bc.
\end{equation}

The set of quantum matrices $A$ are the points of the $r$-dimensional {\em 
quantum algebra} $M_q(r)$, which is defined to be the quotient of the
free algebra in noncommuting variables $x_{ij}$ for $1 \leq i,j, \leq r$,
modulo the left ideal generated by the commutation relations of Equation
\eqref{eq.commutation}.

The algebra $M_q(r)$ has interesting and important structure.
$M_q(r)$ is 
Noetherian, has no zero divisors, and in addition, a basis for the underlying 
vector space is given by the set of {\em sorted monomials} 
$\{\prod_{i,j} a_{ij}^{n_{ij}} | n_{ij} \geq 0\}$
where the product is taken lexicographically; see \cite[Thm IV.4.1]{Ka}. 
An important quotient of $M_q(r)$ is the {\em quantum group} 
$SL_q(r):=M_q(r)/(\det_q-1)$, which is a Hopf algebra \cite[Sec.IV.6]{Ka}
whose representation theory gives rise to the quantum group invariants
of knots, such as the celebrated {\em Jones polynomial}.

Observing that
\begin{eqnarray*}
\tr S^n(A) &=& \sum_{m_1+\dots m_r=n} G(m_1,\dots,m_r) \\
\tr \La^n(A)= &=& \sum_{J \subset \{1,\dots, r\}, |J|=n}
\det_q(A_J)
\end{eqnarray*}
Theorem \ref{thm.1} implies that:

\begin{theorem}
\lbl{thm.2}
If $A$ is in $M_q(r)$, then
$$
\frac{1}{\Ferm(A)}=\sum_{n=0}^\infty \tr S^n(A) 
$$
\end{theorem}

Since the algebra $M_q(r)$ has a vector space basis given by sorted monomials,
it should be possible to give an alternative proof of the quantum MacMahon
Master Theorem using {\em combinatorics on words}, as was done in \cite{FZ} 
for several proofs of the MacMahon Master theorem. We hope to return to this 
alternative point of view in the near future.

\ifx\undefined\bysame
        \newcommand{\bysame}{\leavevmode\hbox
to3em{\hrulefill}\,}
\fi


\begin{thebibliography}{[EMSS]}

\bibitem[A]{A} G. E. Andrews,
        {\em Problems and prospects for basic hypergeometric functions},
        The Theory and Applications of Special Functions (R. Aseky, Editor),
        Academic Press, New York, 1975, 191-224.

\bibitem[FRT]{FRT} L. Fadeev, N. Reshetikhin and L. Takhtadjian,
        {\em Quantization of Lie groups and Lie algebras},
        Leningrad Math. Journal {\bf 1} (1990), 193-225.

\bibitem[FZ]{FZ} D. Foata and D. Zeilberger,
        {\em A combinatorial proof of Bass's evaluation of the Ihara-Selberg
        zeta function for graphs}, 
        Transactions Amer. Math. Soc. {\bf 351} (1999)  2257--2274.

\bibitem[HL]{HL} V. Huynh and TTQ Le,
        {\em On the Colored Jones Polynomial and the Kashaev invariant},
        preprint 2005, {\tt math.GT/0503296}.

\bibitem[Ka]{Ka} C. Kassel,
        {\em Quantum groups}, 
        GTM {\bf 155} Springer-Verlag (1995).

\bibitem[MM]{MM} P. A. MacMahon,
        {\em Combinatory Analysis} vol. 1, Cambridge University Press, 1917.
        Reprinted by Chelsea, 1984.

\bibitem[M]{M} Y. Manin,
        {\em Quantum groups and noncommutative geometry},
        Universit\'e de Montr\'eal, Centre de Recherches Math\'ematiques, 
        Montreal, QC, 1988.

\bibitem[Z1]{Z1} D. Zeilberger, 
        {\em  The algebra of linear partial difference 
        operators and its applications}, 
        SIAM J. Math. Anal. {\bf 11} (1980) 919--934.

\bibitem[Z2]{Z2} D. Zeilberger,
        {\em A holonomic systems approach to special functions identities},
        J. Comput. Appl. Math. {\bf 32} (1990) 321--368.

\end{thebibliography}
\end{document}